\input amstex

\input amstex
\input amsppt.sty
\magnification=\magstep1
 \hsize=36truecc
 \vsize=23.5truecm
\baselineskip=14truept
 \NoBlackBoxes

\def\q{\quad}
\def\qq{\qquad}
\def\mod#1{\ (\text{\rm mod}\ #1)}

\def\t{\text}
\def\qtq#1{\q\t{#1}\q}

\def\f{\frac}
\def\e{\equiv}
\def\a{\alpha}
\def\b{\binom}

\def\o{\omega}

\def\inZ{\in\Bbb Z}

\def\({\left(}
\def\){\right)}
\par\q preprint: December 21, 2009
\par\q
\def\sls#1#2{(\f{#1}{#2})}
 \def\ls#1#2{\big(\f{#1}{#2}\big)}
\def\Ls#1#2{\Big(\f{#1}{#2}\Big)}

\def\ord{\t{\rm ord}}

\let \pro=\proclaim
\let \endpro=\endproclaim

\topmatter
\title Congruences involving $\binom{4k}{2k}$ and $\binom{3k}k$
\endtitle
\author ZHI-Hong Sun\endauthor
\affil School of Mathematical Sciences, Huaiyin Normal University,
\\ Huaian, Jiangsu 223001, PR China
\\ Email: zhihongsun$\@$yahoo.com
\\ Homepage: http://www.hytc.edu.cn/xsjl/szh
\endaffil

 \nologo \NoRunningHeads

\abstract{Let $p$ be a prime greater than $3$. In the paper we
mainly determine $\sum_{k=0}^{[p/4]}\binom{4k}{2k}(-1)^k$,
$\sum_{k=0}^{[p/3]}\binom{3k}k, \sum_{k=0}^{[p/3]}\binom{3k}k(-1)^k$
and $\sum_{k=0}^{[p/3]}\binom{3k}k(-3)^k$ modulo $p$, where $[x]$ is
the greatest integer not exceeding $x$.
\par\q
 }
 \endabstract
 \footnotetext[1] {The author is
supported by the National Natural Sciences Foundation of China
(grant no. 10971078).}
\endtopmatter
\document
\subheading{1. Introduction}
\par\q\par Recently Z.W. Sun and his coauthors studied congruences
involving the binomial coefficient $\b{3k}k$. For example,
Zhao, Pan and Sun [15] showed that for any prime $p>5$,
$$\sum_{k=1}^{p-1}2^k\b{3k}k\e \f 65((-1)^{(p-1)/2}-1)\mod p.$$
In [13] Z.W. Sun investigated $\sum_{k=0}^{p-1}\b{3k}km^{-k}\mod p$
for a prime $p>3$ and $m\not\e 0\mod p$. He gave explicit
congruences in the cases $m=6,7,8,9,13,-\f 14,\f 4{27}, \f 38$.
\par Let $\Bbb Z$ and $\Bbb N$
be the sets of integers and positive integers, respectively. For a
prime $p$ let $\Bbb Z_p$ denote the set of
 those rational numbers whose denominator is not divisible by $p$. Let $p$ be a prime greater than $3$ and $m\in\Bbb Z_p$
with $m\not\e 0\mod p$. In Sections 2 and 3 we study congruences for
$\sum_{k=0}^{[p/4]}\b{4k}{2k}m^{k}$ and
$\sum_{k=0}^{[p/4]}\b{3k}{k}m^{-k}$ modulo $p$ by using Lucas
sequences, binary quadratic forms and the theory of cubic residues,
where $[x]$ is the greatest integer not exceeding $x$. In Section 4
we deduce some congruences involving
$\b{6k}{3k},\b{8k}{4k},\b{10k}{5k},\b {12k}{6k},\b{20k}{10k}$ and
$\b{24k}{12k}$ via formulas for the sum $\sum_{k\e r\mod m}\b nk$ in
the cases $m=3,4,5,6,10,12$. For $a,b,c\in\Bbb Z$ and a prime $p$,
we say that $p$ is represented by $ax^2+bxy+cy^2$ or
$p=ax^2+bxy+cy^2$ if there are integers $x$ and $y$ such that
$p=ax^2+bxy+cy^2$. Let $\sls an$ be the Legendre-Jacobi symbol. As
examples, we have the following typical results.
\par (1.1)  Let $p$ be a prime of the form $4k+3$. Then
$$\sum_{k=0}^{(p-3)/4}(-1)^k\b{4k}{2k}\e
\cases 17^{(p-3)/4}\mod p&\t{if $p\e \pm 1,\pm 4\mod{17}$,}
\\-17^{(p-3)/4}\mod p&\t{if $p\e \pm 2,\pm 8\mod{17}$,}
\\4\cdot 17^{(p-3)/4}\mod p&\t{if $p\e \pm 3,\pm 5\mod{17}$,}
\\-4\cdot 17^{(p-3)/4}\mod p&\t{if $p\e \pm 6,\pm 7\mod{17}$.}
\endcases$$

\par (1.2) Let $p$ be a prime of the form $4k+1$ and so
$p=c^2+d^2$ with $c,d\in\Bbb Z$ and $2\mid d$. Let $a\in\Bbb Z$ with
 $p\nmid (16a^2+1)$. Then
$$\sum_{k=0}^{(p-1)/4}\b{4k}{2k}(-a^2)^k
\e\cases \ls{c-4ad}{16a^2+1}\mod p&\t{if $\ls{16a^2+1}p=1$,}
\\0\mod p&\t{if $\ls{16a^2+1}p=-1$.}
\endcases$$
\par (1.3) Let $p>3$ be a prime. If $\sls p{23}=-1$, then $x\e \sum_{k=0}^{[p/3]}\b{3k}k
 \mod p$ is the unique solution of the congruence $23x^3+3x+1\e
 0\mod p$. If $\sls p{23}=1$, then
$$\sum_{k=0}^{[p/3]}\b{3k}k
 \e\cases 1\mod p&\t{if $p=x^2+xy+52y^2,\ 8x^2+7xy+8y^2$,}
  \\(39x-10y)/(23y)\mod p&\t{if $p=13x^2+xy+4y^2\not=13$,}
 \\-(87x+19y)/(23y)\mod p&\t{if $p=29x^2+5xy+2y^2\not=29$.}
 \endcases$$
 \par (1.4) Let $p$ be a prime with $p\e 13\mod{24}$. Then
$\sum_{k=0}^{(p-1)/12}\b{12k}{6k}\f 1{(-4096)^k}\e 0\mod p.$

\subheading{2. Congruences involving $\b{4k}{2k}$}
\par\q\par For any numbers $P$ and $Q$, let $\{U_n(P,Q)\}$ and $\{V_n(P,Q)\}$ be
the Lucas sequences given by
$$U_0=0,\ U_1=1,\ U_{n+1}=PU_n-QU_{n-1}\ (n\ge 1)$$
and
$$V_0=2,\ V_1=P,\ V_{n+1}=PV_n-QV_{n-1}\ (n\ge 1).$$
It is well known that (see [14])
$$U_n(P,Q)=\cases \f
1{\sqrt{P^2-4Q}}\Big\{\big(\f{P+\sqrt{P^2-4Q}}2\big)^n-\big(\f{P-\sqrt{P^2-4Q}}2\big)^n\Big\}
\\\qq\qq\qq\qq\qq\qq\t{if $P^2-4Q\not=0$,}
\\n\ls P2^{n-1}\qq\qq\qq\qq\t{if $P^2-4Q=0$}
\endcases$$ and
$$V_n(P,Q)=\Ls{P+\sqrt{P^2-4Q}}2^n+\Ls{P-\sqrt{P^2-4Q}}2^n.$$
In particular, we have
$$U_n(2,1)=n\q \t{and}\q U_n(a+b,ab)=\f{a^n-b^n}{a-b} \qtq{for}\ a\neq
b.\tag 2.1$$ As usual, the sequences $F_n=U_n(1,-1)$ and
$L_n=V_n(1,-1)$ are called the Fibonacci sequence and the Lucas
sequence, respectively.  It is easily seen that (see [3, Lemma 1.7])
$$2U_{n+1}(P,Q)=PU_n(P,Q)+V_n(P,Q),\q 2QU_{n-1}(P,Q)=PU_n(P,Q)-V_n(P,Q).\tag 2.2$$

\pro{Lemma 2.1([14, (4.2.39)])} For $n\in\Bbb N$ we have
$$U_{2n+1}(P,Q)=\sum_{k=0}^n\b{n+k}{n-k}(-Q)^{n-k}P^{2k}.$$
\endpro
\pro{Lemma 2.2} Let $p$ be an odd prime and $k\in\{1,2,\ldots, [\f
p4]\}$. Then
$$\b{[\f p4]+k}{[\f p4]-k}\e \b{4k}{2k}\f 1{(-64)^k}\mod p.$$
\endpro
Proof. Suppose $r=1$ or $3$ according as $4\mid p-1$ or $4\mid p-3$.
Then clearly
$$\aligned \b{\f {p-r}4+k}{2k}&
=\f{(\f{p-r}4+k)(\f{p-r}4+k-1)\cdots(\f{p-r}4-k+1)}{(2k)!}
\\&=\f{(p+4k-r)(p+4k-r-4)\cdots(p-(4k+r-4))}{4^{2k}\cdot (2k)!}
\\&\e (-1)^k\f{(4k-r)(4k-r-4)\cdots(4-r)\cdot r(r+4)\cdots(4k+r-4)}
{4^{2k}\cdot (2k)!}
\\&=\f{(-1)^k\cdot (4k)!}{2^{2k}\cdot (2k)!\cdot 4^{2k}\cdot (2k)!}
=\f{\b{4k}{2k}}{(-64)^k} \mod p.
\endaligned$$ So the result follows.

\pro{Lemma 2.3} Let $p$ be an odd prime and
$k\in\{1,2,\ldots,\f{p-1}2\}$. Then
$$\b{(p-1)/2}k\e \f 1{(-4)^k}\b{2k}k\mod {p}.$$
\endpro
Proof. It is clear that
$$\aligned\b{\f{p-1}2}k&=\f{\f{p-1}2(\f{p-1}2-1)\cdots(\f{p-1}2-k+1)}{k!}
=\f{(p-1)(p-3)\cdots(p-(2k-1))}{2^k\cdot k!}
\\&\e \f{(-1)(-3)\cdots(-(2k-1))}{2^k\cdot k!}
=\f{(-1)^k\cdot (2k)!}{(2^k\cdot k!)^2} \mod {p}.\endaligned$$ This
yields the result. \pro{Theorem 2.1} Let $p$ be an odd prime and
$P,Q\in\Bbb Z_p$ with $p\nmid PQ$. Then
$$\sum_{k=0}^{[p/4]}\b{4k}{2k}\Big(\f{P^2}{64Q}\Big)^k\e
(-Q)^{-[p/4]}U_{\f{p+\sls{-1}p}2}(P,Q)\mod p$$ and
$$\sum_{k=0}^{[p/4]}\b{4k}{2k}\Ls Q{4P^2}^k\e \Ls{P}p
U_{\f{p+1}2}(P,Q)\mod p.$$
\endpro
Proof. Using
 Lemmas 2.1 and 2.2 we see that
$$\aligned U_{2[\f
p4]+1}(P,Q)&=\sum_{k=0}^{[p/4]}\b{[\f p4]+k}{[\f p4]-k}(-Q)^{[\f
p4]-k}P^{2k}
\\&\e(-Q)^{[\f
p4]}\sum_{k=0}^{[\f p4]}\f{\b{4k}{2k}}{(-64)^k}\Ls{P^2}{-Q}^k\mod
p.\endaligned$$ Note that $2[\f p4]+1=(p+\sls{-1}p)/2$. We deuce the
first result.
\par Using Lemma 2.3 we see that
$$\aligned V_{\f{p-1}2}(P,(P^2-4Q)/4)&=\Ls{P+2\sqrt Q}2^{\f{p-1}2}
+\Ls{P-2\sqrt Q}2^{\f{p-1}2}
\\&=\f
1{2^{\f{p-1}2}}\sum_{k=0}^{\f{p-1}2}\b{\f{p-1}2}kP^{\f{p-1}2-k}\big((2\sqrt
Q)^k+(-2\sqrt Q)^k\big)
\\&=\f
2{2^{\f{p-1}2}}\sum_{k=0}^{[p/4]}\b{\f{p-1}2}{2k}P^{\f{p-1}2-2k}(2\sqrt
Q)^{2k}
\\&\e \f
2{2^{\f{p-1}2}}\sum_{k=0}^{[p/4]}\b{4k}{2k}(-4)^{-2k}P^{\f{p-1}2-2k}(4Q)^k
\\&\e 2\Ls {2P}p\sum_{k=0}^{[p/4]}\b{4k}{2k}\Ls Q{4P^2}^k\mod
p.\endaligned$$ By appealing to [7, Lemma 3.1] we have
$$U_{\f{p+1}2}(P,Q)\e \f 12\Ls 2pV_{\f{p-1}2}(P,(P^2-4Q)/4)
\e \Ls Pp\sum_{k=0}^{[p/4]}\b{4k}{2k}\Ls Q{4P^2}^k\mod p.$$ This
completes the proof.
 \pro{Lemma 2.4} For $n\in\Bbb N$ we have
$U_n(1,1)=(-1)^{n-1}\sls n3.$\endpro
Proof. Set $\omega
=(-1+\sqrt{-3})/2.$ Then
$$\aligned U_n(1,1)&=\f
1{\sqrt{-3}}\Big\{\ls{1+\sqrt{-3}}2^n-\ls{1-\sqrt{-3}}2^n\Big\}
\\&=\f
{(-1)^n}{\sqrt{-3}}\Big\{\ls{-1-\sqrt{-3}}2^n-\ls{-1+\sqrt{-3}}2^n\Big\}
\\&=\f{(-1)^n}{\omega (1-\omega )}\(\omega^{2n}-\omega^n\)
=\cases 0&\t{if}\ 3\mid n
\\(-1)^{n-1}&\t{if}\ 3\mid n-1\\(-1)^n&\t{if}\ 3\mid
n-2.\endcases\endaligned$$ This yields the result.

 \pro{Corollary 2.1} Let $p$ be an odd prime.
Then
$$\aligned &\sum_{k=0}^{[p/4]}\b{4k}{2k}\f 1{16^k}\e\f 12\Ls 2p\mod
p,\\&\sum_{k=0}^{[p/4]}\b{4k}{2k}\f 1{4^k}\e \cases
(-1)^{\f{p-1}2}\mod p&\t{if $p\e 1\mod 3$,}
\\0\mod p&\t{if $p\e 2\mod 3$}\endcases\endaligned$$ and
$$\sum_{k=0}^{[p/4]}\b{4k}{2k}\f 1{64^k}\e\cases 0\mod p&\t{if $p\e
5,7,17,19\mod{24}$,}
\\1\mod p&\t{if $p\e 1,23\mod{24}$,}
\\-1\mod p&\t{if $p\e 11,13\mod{24}$.}\endcases$$
\endpro
Proof. Taking $(P,Q)=(2,1),(1,1)$ in Theorem 2.1 and then applying
(2.1) and Lemma 2.4 we deduce the result.
 \pro{Theorem 2.2} Let $p$ be an odd
prime and $x\in\Bbb Z_p$ with $x\not\e 0,1\mod p$.
\par $(\t{\rm i})$ If $p\e 1\mod 4$, then
$$\sum_{k=0}^{[p/4]}\b{4k}{2k}\Ls x{16}^k\e x^{\f{p-1}4}\sum_{k=0}^{[p/4]}\b{4k}{2k}
\f 1{(16x)^k}\mod p.$$

\par $(\t{\rm ii})$ If $p\e 3\mod 4$, then
$$\sum_{k=0}^{[p/4]}\b{4k}{2k}\f 1{(16x)^k}\e \Big(1-\f 1x\Big)^{\f{p-3}4}\sum_{k=0}^{[p/4]}\b{4k}{2k}
\f 1{(16(1-x))^k}\mod p.$$
\endpro
Proof. If $p\e 1\mod 4$, by Theorem 2.1 we have
$$\sum_{k=0}^{[p/4]}\b{4k}{2k}\Ls x{16}^k\e\Ls 2pU_{\f{p+1}2}(2,x)\e
(-1)^{\f{p-1}4}\cdot (-x)^{\f{p-1}4}\sum_{k=0}^{[p/4]}\b{4k}{2k} \f
1{(16x)^k}\mod p.$$ If $p\e 3\mod 4$, from [7, Lemma 3.1] we know
that $U_{\f{p-1}2}(2,x)\e -\sls 2pU_{\f{p-1}2}(2,1-x)\mod p.$
 Now applying Theorem 2.1 and the fact $\sls 2p=-(-1)^{(p-3)/4}$ we deduce
$$\align\sum_{k=0}^{[p/4]}\b{4k}{2k} \f
1{(16x)^k}&\e (-x)^{-[p/4]}U_{\f{p-1}2}(2,x)\e
x^{-[p/4]}U_{\f{p-1}2}(2,1-x)
\\&\e (x-1)^{[\f p4]}x^{-[\f p4]}
\sum_{k=0}^{[p/4]}\b{4k}{2k} \f 1{(16(1-x))^k}\mod p.\endalign$$ So
the theorem is proved.
 \pro{Corollary 2.2} Let $p$ be a prime of the form $8k+7$. Then
$$\sum_{k=0}^{(p-3)/4}\b{4k}{2k}8^{-k}\e 0\mod p.$$
\endpro
Proof. Taking $x=\f 12$ in Theorem 2.2(ii) we deduce the result.
 \pro{Theorem 2.3} Let
$p$ be an odd prime and $P,Q\in\Bbb Z_p$ with $p\nmid PQ(P^2-4Q)$
and $\sls Qp=1$. \par $(\t{\rm i})$ If $\sls{4Q-P^2}p=-1$, then
$\sum_{k=0}^{[p/4]}\b{4k}{2k}\sls{P^2}{64Q}^k\e 0\mod p.$
 \par $(\t{\rm ii})$ If $\sls{P^2-4Q}p=-1$, then
$\sum_{k=0}^{[p/4]}\b{4k}{2k}\sls Q{4P^2}^k\e 0\mod p.$
\endpro
Proof. Since $\sls Qp=1$, it is well known that (see [3])
$U_{(p-\sls{P^2-4Q}p)/2}(P,Q)\e 0\mod p$. This together with Theorem
2.1 yields the result.
 \pro{Theorem 2.4} Let
$p>5$ be a prime. Then
$$\sum_{k=0}^{[p/4]}\f{\b{4k}{2k}}{(-16)^k}
\e \cases (-1)^{\f{p-1}8}2^{\f{p-1}4}\mod p&\t{if $p\e 1\mod{8}$,}
\\(-1)^{\f{p-3}8}2^{\f{p-3}4}\mod p&\t{if $p\e
3\mod{8}$,}
\\0\mod p&\t{if $p\e
5\mod{8}$,}
\\(-1)^{\f{p+1}8}2^{\f{p-3}4}\mod p&\t{if $p\e
7\mod{8}$.}\endcases$$
\endpro
Proof. Taking $P=2$ and $Q=-1$ in Theorem 2.1 we obtain
$$\sum_{k=0}^{[p/4]}\b{4k}{2k}(-16)^{-k}\e
U_{(p+\sls{-1}p)/2}(2,-1)\mod p.$$ Now applying [4, Theorem 2.3] we
deduce the result.
 \pro{Theorem 2.5} Let $p>5$
be a prime. Then
$$\sum_{k=0}^{[p/4]}\f{\b{4k}{2k}}{(-64)^k}
\e \cases (-1)^{[\f{p+5}{10}]}5^{[\f p4]}\mod p&\t{if $p\e
1,3,7,9\mod{20}$,}
\\2(-1)^{[\f{p+5}{10}]}5^{\f{p-3}4}\mod p&\t{if $p\e
11,19\mod{20}$,}
\\0\mod p&\t{if $p\e
13,17\mod{20}$}\endcases$$ and
$$\sum_{k=0}^{[p/4]}\f{\b{4k}{2k}}{(-4)^k}
\e \cases (-1)^{[\f{p+5}{10}]}5^{[\f p4]}\mod p&\t{if $p\e
1,9,11,19\mod{20}$,}
\\-2(-1)^{[\f{p+5}{10}]}5^{\f{p-3}4}\mod p&\t{if $p\e
3,7\mod{20}$,}
\\0\mod p&\t{if $p\e
13,17\mod{20}$.}\endcases$$

\endpro
Proof. Taking $P=1$ and $Q=-1$ in Theorem 2.1 we obtain
$$\sum_{k=0}^{[p/4]}\b{4k}{2k}(-64)^{-k}\e
F_{\f{p+\sls{-1}p}2}\mod
p\qtq{and}\sum_{k=0}^{[p/4]}\b{4k}{2k}(-4)^{-k}\e F_{\f{p+1}2}\mod
p.$$ Now applying [11, Corollaries 1 and 2] we deduce the result.

\pro{Theorem 2.6} Let $p$ be an odd prime with $p\not =17$.
\par $(\t{\rm i})$ If $p\e 1\mod 4$, then
$$\sum_{k=0}^{[p/4]}(-1)^k\b{4k}{2k}\e
\cases 0\mod p&\t{if $p\e \pm 3,\pm 5,\pm 6,\pm 7\mod{17}$,}
\\17^{(p-1)/4}\mod p&\t{if $p\e \pm 1,\pm 4\mod{17}$,}
\\-17^{(p-1)/4}\mod p&\t{if $p\e \pm 2,\pm 8\mod{17}$.}
\endcases$$
\par $(\t{\rm ii})$ If $p\e 3\mod 4$, then
$$\sum_{k=0}^{[p/4]}(-1)^k\b{4k}{2k}\e
\cases 17^{(p-3)/4}\mod p&\t{if $p\e \pm 1,\pm 4\mod{17}$,}
\\-17^{(p-3)/4}\mod p&\t{if $p\e \pm 2,\pm 8\mod{17}$,}
\\4\cdot 17^{(p-3)/4}\mod p&\t{if $p\e \pm 3,\pm 5\mod{17}$,}
\\-4\cdot 17^{(p-3)/4}\mod p&\t{if $p\e \pm 6,\pm 7\mod{17}$.}
\endcases$$
\endpro
Proof. Taking $P=8$ and $Q=-1$ in Theorem 2.1 we see that
$$\sum_{k=0}^{[p/4]}(-1)^k\b{4k}{2k}\e U_{\f{p+\sls
{-1}p}2}(8,-1)\mod p.\tag 2.3$$ By (2.2) we have
$U_{\f{p+1}2}(8,-1)=4U_{\f{p-1}2}(8,-1)+\f 12V_{\f{p-1}2}(8,-1)$.
From the above and [10, Corollary 4.5] we deduce the result.
\pro{Theorem 2.7} Let $p$ be an odd prime with $p\not =3,13$. Then
$$\aligned&\Ls 3p\sum_{k=0}^{[p/4]}\b{4k}{2k}\f 1{(-36)^k}
\\&\e \cases  13^{\f{p-1}4}\mod p&\t{if $p\e
1,9,-23\mod{52}$,}\\
-13^{\f{p-1}4}\mod p&\t{if $p\e
-3,17,25\mod{52}$,}\\
3\cdot 13^{\f{p-3}4}\mod p&\t{if $p\e
-1,-9,23\mod{52}$,}\\
-3\cdot 13^{\f{p-3}4}\mod p&\t{if $p\e
3,-17,-25\mod{52}$,}\\
0\mod p&\t{if $p\e
5,-7,-11,-15,-19,21\mod{52}$,}\\
2\cdot 13^{\f{p-3}4}\mod p&\t{if $p\e
-5,7,11\mod{52}$,}\\
-2\cdot 13^{\f{p-3}4}\mod p&\t{if $p\e
15,19,-21\mod{52}$.}\endcases\endaligned$$
\endpro
Proof. Taking $P=3$ and $Q=-1$ in Theorem 2.1 we see that
$$\sum_{k=0}^{[p/4]}\b{4k}{2k}\f 1{(-36)^k}\e \Ls 3pU_{\f{p+1}2}(3,-1)\mod p.$$
Now applying [7, Corollary 4.1] we deduce the result.

\pro{Theorem 2.8} Let $p$ be a prime such that $p\e
1,9,11,13,19,37\mod{40}$. Then
$$\aligned&\sum_{k=0}^{[p/4]}\b{4k}{2k}\f 1{(-144)^k}
\\&\e
\cases (-1)^{[\f p3]+\f y2}\mod p&\t{if $p=x^2+10y^2\e
1,9\mod{40}$,}
\\(-1)^{[\f p3]}\f yx\mod p&\t{if $p=x^2+10y^2\e
11,19\mod{40}$ and $4\mid x-y$,}
\\(-1)^{[\f p3]+\f y2}\mod p&\t{if $p=5x^2+2y^2\e
13,37\mod{40}$.}
\endcases\endaligned$$
\endpro
Proof. Taking $P=6$ and $Q=-1$ in Theorem 2.1 we see that
$$\sum_{k=0}^{[p/4]}\b{4k}{2k}\f 1{(-144)^k}\e \Ls 6pU_{\f{p+1}2}(6,-1)\mod p.$$
Now applying [10, Theorem 5.4] we deduce the result.

 \pro{Lemma 2.5 ([5, Lemma 3.4])} Let $p$ be an odd prime
and $P,Q\in\Bbb Z_p$ with $p\nmid Q(P^2-4Q)$. If $\ls Qp=1$ and
$c^2\e Q\mod p$ for $c\in\Bbb Z_p$, then
$$U_{\f{p+1}2}(P,Q)\e\cases \ls{P-2c}p\mod p&\t{if
$\sls{P^2-4Q}p=1$,}
\\0\mod p&\t{if
$\sls{P^2-4Q}p=-1$}
\endcases$$ and
$$U_{\f{p-1}2}(P,Q)\e\cases 0\mod p&\t{if
$\sls{P^2-4Q}p=1$,}
\\\f 1c\ls{P-2c}p\mod p&\t{if
$\sls{P^2-4Q}p=-1$.}
\endcases$$
\endpro
\pro{Theorem 2.9} Let $p$ be an odd  prime and $a\in\Bbb Z_p$ with
$p\nmid (16a^2-1)$. Then
$$\sum_{k=0}^{[p/4]}\b{4k}{2k}a^{2k}\e
\cases 0\mod p&\t{if $\sls{1-16a^2}p=-1$,}
\\ \ls{1-4a}p\mod p&\t{if $\sls{1-16a^2}p=1$.}
\endcases$$
\endpro
Proof. Putting $P=8a$ and $Q=1$ in Theorem 2.1 we deduce that
$$\sum_{k=0}^{[p/4]}\b{4k}{2k}a^{2k}\e
(-1)^{[p/4]}U_{\f{p+\sls{-1}p}2}(8a,1)\mod p.$$ If $p\e 1\mod 4$, by
Lemma 2.5 we have
$$(-1)^{\f{p-1}4}U_{\f{p+1}2}(8a,1)
\e\cases (-1)^{\f{p-1}4}\ls {8a-2}p=\ls {1-4a}p\mod p&\t{if
$\sls{16a^2-1}p=1$,}\\ 0\mod p&\t{if $\sls{16a^2-1}p=-1$.}
\endcases$$
If $p\e 3\mod 4$, by Lemma 2.5 we have
$$(-1)^{\f{p-3}4}U_{\f{p-1}2}(8a,1)
\e\cases (-1)^{\f{p-3}4}\ls {8a-2}p=\ls {1-4a}p\mod p&\t{if
$\sls{16a^2-1}p=-1$,}\\ 0\mod p&\t{if $\sls{16a^2-1}p=1$.}
\endcases$$
Now combining all the above we obtain the result.

\pro{Corollary 2.2} Let $p>5$ be a prime. Then
$$\sum_{k=0}^{[p/4]}\b{4k}{2k}\e
\cases 0\mod p&\t{if $p\e 7,11,13,14\mod{15}$,}
\\ 1\mod p&\t{if $p\e 1,4\mod{15}$,}
\\ -1\mod p&\t{if $p\e 2,8\mod{15}$.}
\endcases$$
\endpro
Proof. Since $\sls{-15}p=\sls p{15}=\sls p3\sls p5$, putting $a=1$
in Theorem 2.9 we deduce the result.
 \pro{Corollary 2.3} Let $p>7$
be a prime. Then
$$\sum_{k=0}^{[p/4]}\b{4k}{2k}2^{2k}\e
\cases 1\mod p&\t{if $p\e 1,2,4\mod 7$,}
\\ 0\mod p&\t{if $p\e 3,5,6\mod 7$.}
\endcases$$
\endpro
Proof. Since $\sls{-63}p=\sls {-7}p=\sls p7$, putting $a=2$ in
Theorem 2.9 we deduce the result.

\pro{Corollary 2.4} Let $p>3$ be a prime. Then
$$\sum_{k=0}^{[p/4]}\b{4k}{2k}\f 1{2^{2k}}\e
\cases (-1)^{(p-1)/2}\mod p&\t{if $p\e 1\mod 3$,}
\\ 0\mod p&\t{if $p\e 2\mod 3$.}
\endcases$$
\endpro
Proof. Taking $a=\f 12$ in Theorem 2.9 we deduce the result.

\pro{Theorem 2.10} Let $p$ be a prime of the form $4k+1$ and
$p=c^2+d^2$ with $c,d\in\Bbb Z$ and $2\mid d$. Let $b,m\in\Bbb Z$
with $\t{gcd}(b,m)=1$ and $p\nmid m(b^2+4m^2)$. Then
$$\aligned&\Ls mp\sum_{k=0}^{(p-1)/4}\b{4k}{2k}\Big(-\f{b^2}{64m^2}\Big)^k
\e \Ls bp\sum_{k=0}^{(p-1)/4}\b{4k}{2k}\Big(-\f{m^2}{4b^2}\Big)^k
\\&\e\cases \ls{bc+2md}{b^2+4m^2}\mod p&\t{if $2\nmid b$ and
$\sls{b^2+4m^2}p=1$,}
\\(-1)^{\f{(\f b2c+md)^2-1}8+\f d2}\ls{\f b2c+md}{(\sls
b2^2+m^2)/2}\mod p&\t{if $2\ \Vert\ b$ and $\sls{b^2+4m^2}p=1$,}
\\\ls{mc-\f b2d}{\f{b^2}4+m^2}\mod p
&\t{if $4\mid b$ and $\sls{b^2+4m^2}p=1$,}
\\0\mod p&\t{if $\sls{b^2+4m^2}p=-1$.}
\endcases\endaligned$$
In particular, taking $b=8a$ and $m=1$ we have
$$\aligned\sum_{k=0}^{(p-1)/4}\b{4k}{2k}(-a^2)^k&\e \Ls {2a}p\sum_{k=0}^{(p-1)/4}\b{4k}{2k}
\f 1{(-256a^2)^k} \\&\e\cases \ls{c-4ad}{16a^2+1}\mod p&\t{if
$\ls{16a^2+1}p=1$,}
\\0\mod p&\t{if $\ls{16a^2+1}p=-1$.}
\endcases\endaligned$$
\endpro
Proof. Putting $P=b$ and $Q=-m^2$ in Theorem 2.1 we see that
$$U_{\f{p+1}2}(b,-m^2)\e m^{\f{p-1}2}
\sum_{k=0}^{\f{p-1}4}\b{4k}{2k}\Big(-\f{b^2}{64m^2}\Big)^k \e\Ls
bp\sum_{k=0}^{\f{p-1}4}\b{4k}{2k}\Big(-\f{m^2}{4b^2}\Big)^k\mod p.$$
Now applying [10, Theorem 3.2] we deduce the result.

\pro{Corollary 2.5} Let $p$ be a prime with $p\e 1\mod 4$ and $p\e
\pm 1,\pm 2,\pm 4,\pm 8$ $\mod{17}$. Let $p=c^2+d^2$ with
$c,d\in\Bbb Z$ and $2\mid d$. Then
$$\sum_{k=0}^{(p-1)/4}(-1)^k\b{4k}{2k}\e\Ls{c-4d}{17}\mod p.$$
\endpro
Proof. Taking $a=1$ in Theorem 2.10 we obtain the result.

\pro{Theorem 2.11} Let $p$ be an odd prime, and $b,m\in\Bbb Z$ with
$p\nmid bm(b^2+4m^2)$. Then there is a unique $\delta_p\in\{1,-1\}$
such that
$$\aligned&\Ls
bp\sum_{k=0}^{[p/4]}\b{4k}{2k}\Big(-\f{m^2}{4b^2}\Big)^k
\\&\e \cases
\delta_p(b^2+4m^2)^{\f{p-1}4}\mod p&\t{if $4\mid p-1$ and
$\sls{b^2+4m^2}p=1$,}
\\b\delta_p(b^2+4m^2)^{\f{p-3}4}\mod p&\t{if $4\mid p-3$ and
$\sls{b^2+4m^2}p=1$,}\\0\mod p  &\t{if $4\mid p-1$ and
$\sls{b^2+4m^2}p=-1$,}
\\2m\delta_p(b^2+4m^2)^{\f{p-3}4}\mod p&\t{if $4\mid p-3$ and
$\sls{b^2+4m^2}p=-1$}\endcases\endaligned$$ and
$$\aligned&\Ls
mp\sum_{k=0}^{[p/4]}\b{4k}{2k}\Big(-\f{b^2}{64m^2}\Big)^k
\\&\e \cases
\delta_p(b^2+4m^2)^{\f{p-1}4}\mod p&\t{if $4\mid p-1$ and
$\sls{b^2+4m^2}p=1$,}
\\2m\delta_p(b^2+4m^2)^{\f{p-3}4}\mod p&\t{if $4\mid p-3$ and
$\sls{b^2+4m^2}p=1$,}\\0\mod p  &\t{if $4\mid p-1$ and
$\sls{b^2+4m^2}p=-1$,}
\\-b\delta_p(b^2+4m^2)^{\f{p-3}4}\mod p&\t{if $4\mid p-3$ and
$\sls{b^2+4m^2}p=-1$.}\endcases\endaligned$$ Moreover, if $p'$ is
also an odd prime satisfying $p'\nmid bm(b^2+4m^2)$ and $p'\e \pm
p\mod{(3-(-1)^b)(b^2+4m^2)}$, then $\delta_p=\delta_{p'}$. Indeed,
$$\delta_p=\cases \ls{b+2mi}p_4&\t{if $\ls{b^2+4m^2}p=1$,}
\\\ls{b+2mi}p_4i&\t{if $\ls{b^2+4m^2}p=-1$,}
\endcases$$
where $\sls{b+ci}p_4$ is the quartic Jacobi symbol.
\endpro
Proof. Putting $P=b$ and $Q=-m^2$ in Theorem 2.1 and then applying
[7, Theorem 4.1] we deduce the result.

\pro{Theorem 2.12} Let $p$ be a prime of the form $4k+1$ and
$a\in\Bbb Z$ with $p\nmid (1+16a^2)(1-16a^2)$. Let $p=c^2+d^2$ with
$c,d\in\Bbb Z$ and $2\mid d$. Then
$$2\sum_{k=0}^{[p/8]}\b{8k}{4k}a^{4k}\e
\cases  \sls{1-4a}p+\sls{c-4ad}{16a^2+1}\mod p&\t{if
$\sls{1-16a^2}p=\sls{1+16a^2}p=1$,}
\\\sls {1-4a}p\mod p&\t{if
$\sls{1-16a^2}p=-\sls{1+16a^2}p=1$,}
\\\sls {c-4ad}{16a^2+1}\mod p&\t{if
$\sls{1-16a^2}p=-\sls{1+16a^2}p=-1$,}
\\0\mod p&\t{if $\sls{1-16a^2}p=\sls{1+16a^2}p=-1$.}
\endcases$$
\endpro
Proof. Since
$$2\sum_{k=0}^{[p/8]}\b{8k}{4k}a^{2k}=\sum_{k=0}^{(p-1)/4}\b{4k}{2k}a^{2k}+
\sum_{k=0}^{(p-1)/4}\b{4k}{2k}(-1)^ka^{2k},$$ from Theorems 2.9 and
2.10 we deduce the result.
 \pro{Corollary 2.7} Let $p$ be a prime of the form
$4k+1$ with $p\not=5,17$. Let $p=c^2+d^2$ with $c,d\in\Bbb Z$ and
$2\mid d$. Then
$$2\sum_{k=0}^{[p/8]}\b{8k}{4k}\e
\cases  \sls p3+\sls{c-4d}{17}\mod p&\t{if $\sls p{15}=\sls
p{17}=1$,}
\\\sls p3\mod p&\t{if
$\sls p{15}=-\sls p{17}=1$,}
\\\sls {c-4d}{17}\mod p&\t{if
$\sls p{15}=-\sls p{17}=-1$,}
\\0\mod p&\t{if $\sls p{15}=\sls p{17}=-1$.}
\endcases$$
\endpro
Proof. Taking $a=1$ in Theorem 2.12 we obtain the result.
\subheading{3. Congruences involving $\b{3k}k$}
 \pro{Lemma 3.1} Let
$p>3$ be a prime and $k\in\{1,2,\ldots, [\f p3]\}$. Then
$$\b{[\f p3]+k}{[\f p3]-k}\e \b{3k} k\f 1{(-27)^k}\mod p.$$
\endpro
Proof. Suppose $r=1$ or $2$ according as $3\mid p-1$ or $3\mid p-2$.
Then clearly
$$\aligned \b{\f {p-r}3+k}{2k}&
=\f{(\f{p-r}3+k)(\f{p-r}3+k-1)\cdots(\f{p-r}3-k+1)}{(2k)!}
\\&=\f{(p+3k-r)(p+3k-r-3)\cdots(p-(3k+r-3))}{3^{2k}\cdot (2k)!}
\\&\e (-1)^k\f{(3k-r)(3k-r-3)\cdots(3-r)\cdot r(r+3)\cdots(3k+r-3)}
{3^{2k}\cdot (2k)!}
\\&=\f{(-1)^k\cdot (3k)!}{3\cdot 6\cdots 3k\cdot 3^{2k}\cdot (2k)!}
=\f{(-1)^k\cdot(3k)!}{3^k\cdot k!\cdot 3^{2k}\cdot (2k)!} \mod p.
\endaligned$$ So the result follows.
 \pro{Theorem 3.1} Let $p>3$ be a prime
and $a,b\in\Bbb Z_p$ with $p\nmid ab$. Then
$$\sum_{k=0}^{[\f
p3]}\b{3k}k\f{b^{2k}}{a^k}\e (-3a)^{-[p/3]}U_{2[\f p3]+1}(9b,3a)\mod
p.$$\endpro Proof. Using Lemmas 2.1 and 3.1 we see that for
$P,Q\in\Bbb Z_p$ with $p\nmid PQ$,
$$\aligned U_{2[\f
p3]+1}(P,Q)&=\sum_{k=0}^{[p/3]}\b{\f{p-r}3+k}{\f{p-r}3-k}(-Q)^{\f{p-r}3-k}P^{2k}
\\&\e(-Q)^{[\f
p3]}\sum_{k=0}^{[\f p3]}\b{3k}k\f 1{27^k}\Ls{P^2}Q^k\mod
p.\endaligned\tag 3.1$$ Now taking $P=9b$ and $Q=3a$ in (3.1) we
deduce the result.

\pro{Theorem 3.2} Let $p>3$ be a prime. Then

$$\sum_{k=0}^{[p/3]}\f{\b{3k}k}{27^k}\e\cases
1\mod p&\t{if}\ p\e \pm 1\mod 9,\\-1\mod p&\t{if}\ p\e \pm 2\mod
9,\\0\mod p &\t{if}\ p\e \pm 4\mod 9.\endcases
$$
\endpro
Proof.
 Taking $a=\f 13$ and $b=\f 19$ in Theorem 3.1 and then applying
Lemma 2.3 we deduce $$\sum_{k=0}^{[p/3]}\b{3k}k\f 1{27^k}\e
(-1)^{[\f p3]}\Ls{2[\f p3]+1}3\mod p.$$ This yields the result.

 \pro{Lemma 3.2} Let $p>3$ be a
prime and $P,Q\in\Bbb Z_p$ with $p\nmid PQ$. Then

$$U_{2[\f p3]+1}(P,Q)\e\cases -Q^{1-\f{p-\sls p3}3}U_{\f{p-\sls p3}3-1}(P,Q)\mod p&\t{if}\ \sls{P^2-4Q}p=1,
\\-Q^{-\f{p-\sls p3}3}U_{\f{p-\sls p3}3+1}(P,Q)\mod p&\t{if}\ \sls{P^2-4Q}p=-1.\endcases$$\endpro
Proof. Since $2[\f p3]+1=p-\f{p-\sls p3}3$ and
$\f{P\pm\sqrt{P^2-4Q}}2=\f Q{(P\mp\sqrt{P^2-4Q})/2}$, we see that

$$\aligned U_{2[\f p3]+1}(P,Q)&=\f 1{\sqrt{P^2-4Q}}\Big\{\Ls{P+\sqrt{P^2-4Q}}2^{2[\f p3]+1}
- \Ls{P-\sqrt{P^2-4Q}}2^{2[\f p3]+1}\Big\}
\\&=\f
1{\sqrt{P^2-4Q}}\Big\{\Ls{P+\sqrt{P^2-4Q}}2^p\Ls{(P-\sqrt{P^2-4Q})/2}Q^{\f{p-\sls
p3}3}\\&\qq-
\Ls{P-\sqrt{P^2-4Q}}2^p\Ls{(P+\sqrt{P^2-4Q})/2}Q^{\f{p-\sls
p3}3}\Big\}.\endaligned$$ Since
$$\aligned\Ls{P\pm\sqrt{P^2-4Q}}2^p&\e\f{P^p\pm(\sqrt{P^2-4Q})^p}{2^p}
\e\f{P\pm\sqrt{P^2-4Q}(P^2-4Q)^{\f{p-1}2}}2
\\&\e\f{P\pm\sls{P^2-4Q}p\sqrt{P^2-4Q}}2\mod p,\endaligned$$ by the above we have
$$\aligned U_{2[\f p3]+1}(P,Q)&\e\f{Q^{-(p-\sls p3)/3}}{\sqrt{P^2-4Q}}
\Big\{\f{P+\sls{P^2-4Q}p\sqrt{P^2-4Q}}2\Ls{P-\sqrt{P^2-4Q}}2^{\f{p-\sls
p3}3}
\\&\qq-\f{P-\sls{P^2-4Q}p\sqrt{P^2-4Q}}2\Ls{P+\sqrt{P^2-4Q}}2^{\f{p-\sls
p3}3}\Big\}\mod p.
\endaligned$$
If $\sls{P^2-4Q}p=-1$, by the above we have
$$U_{2[\f p3]+1}(P,Q)\e -Q^{-\f{p-\sls p3}3}U_{\f{p-\sls p3}3+1}(P,Q)\mod p.$$ If
$\sls{P^2-4Q}p=1$, by the above and the fact $\f{P\pm
\sqrt{P^2-4Q}}2=\f Q{(P\mp\sqrt{P^2-4Q})/2}$ we have
$$U_{2[\f p3]+1}(P,Q)\e -Q^{1-(p-\sls p3)/3}U_{\f{p-\sls p3}3-1}(P,Q)\mod p.$$ So the lemma is proved.

\pro{Theorem 3.3} Let $p>3$ be a prime and $a,b\in\Bbb Z_p$ with
$p\nmid ab$. Then
$$\sum_{k=0}^{[p/3]}\b{3k}k\f{b^{2k}}{a^k}\e\cases
(-3a)^{[\f p3]+1}U_{\f{p-\sls p3}3-1}(9b,3a)\mod p&\t{if}\
\sls{81b^2-12a}p=1,\\
-(-3a)^{[\f p3]}U_{\f{p-\sls p3}3+1}(9b,3a)\mod p&\t{if}\
\sls{81b^2-12a}p=-1.\endcases$$\endpro Proof. From Theorem 3.1 and
Lemma 3.2 we deduce
$$\aligned &\sum_{k=0}^{[p/3]}\b{3k}k\f{b^{2k}}{a^k}\e (-3a)^{-[\f p3]}U_{2[\f p3]+1}(9b,3a)
\\&\e\cases
-(-3a)^{-[\f p3]}\cdot (3a)^{1-\f{p-\sls p3}3}U_{\f{p-\sls
p3}3-1}(9b,3a)\mod p&\t{if}\ \Ls{81b^2-12a}p=1,\\ -(-3a)^{-[\f
p3]}\cdot (3a)^{-\f{p-\sls p3}3}U_{\f{p-\sls p3}3+1}(9b,3a)\mod
p&\t{if}\ \Ls{81b^2-12a}p=-1.\endcases\endaligned$$ To see the
result we note that $2[\f p3]=p-1-\f{p-\sls p3}3$ and so
$$(3a)^{-[\f p3]-\f{p-\sls p3}3}=(3a)^{[\f p3]-(p-1)}\e (3a)^{[\f
p3]}\mod p.$$ \pro{Corollary 3.1} Let $p>5$ be a prime. Then
$$\sum_{k=0}^{[p/3]}\f{\b{3k}k}{(-27)^k}\e
\cases F_{\f{p-\sls p3}3-1}\mod p&\t{if $\sls p5=1$,}
\\-F_{\f{p-\sls p3}3+1}\mod p&\t{if $\sls p5=-1$.}
\endcases$$\endpro
Proof. Taking $a=-\f 13$ and $b=\f 19$ in Theorem 3.3 we obtain the
result.
 \pro{Theorem 3.4} Let $p>5$ be a prime, and let $\varepsilon_p=1,-1,0$ according as $p\e \pm 1,\pm 2$
 or $\pm 4\mod 9$.
\par $(\t{\rm i})$ If $p\e 1,4\mod{15}$ and so $p=x^2+15y^2$ with
$x,y\in\Bbb Z$, then
$$2\sum_{k=0}^{[p/6]}\f{\b{6k}{2k}}{27^{2k}}-\varepsilon_p
\e \sum_{k=0}^{[p/3]}\f{\b{3k}k}{(-27)^k}\e \cases 1\mod p&\t{if
$3\mid y$,}
\\(x-5y)/(10y)\mod p&\t{if $3\mid y-x$.}
\endcases$$
\par $(\t{\rm ii})$ If $p\e 2,8\mod{15}$ and so $p=5x^2+3y^2$ with
$x,y\in\Bbb Z$, then
$$2\sum_{k=0}^{[p/6]}\f{\b{6k}{2k}}{27^{2k}}-\varepsilon_p
\e\sum_{k=0}^{[p/3]}\f{\b{3k}k}{(-27)^k}\e \cases 1\mod p&\t{if
$3\mid y$,}
\\-(x+y)/(2y)\mod p&\t{if $3\mid y-x$.}
\endcases$$
\endpro
Proof. By Theorem 3.2 we have
$$2\sum_{k=0}^{[p/6]}\f{\b{6k}{2k}}{27^{2k}}=\sum_{k=0}^{[p/3]}\Big(
\f{\b{3k}k}{27^k}+\f{\b{3k}k}{(-27)^k}\Big)\e
\varepsilon_p+\sum_{k=0}^{[p/3]}\f{\b{3k}k}{(-27)^k}\mod p.$$ If
$p=x^2+15y^2\e 1,4\mod{15}$, by [6, Theorem 6.2] we have
$$F_{\f{p-1}3}\e \cases 0\mod p&\t{if $3\mid y$,}
\\-x/(5y)\mod p&\t{if $3\mid y-x$}
\endcases\ \t{and}\
L_{\f{p-1}3}\e \cases 2\mod p&\t{if $3\mid y$,}
\\-1\mod p&\t{if $3\nmid y$.}
\endcases\tag 3.2$$
If $p=5x^2+3y^2\e 2,8\mod{15}$, by [6, Theorem 6.2] we have
$$F_{\f{p+1}3}\e \cases 0\mod p&\t{if $3\mid y$,}
\\ x/y\mod p&\t{if $3\mid y-x$}
\endcases\ \t{and}\
L_{\f{p+1}3}\e \cases -2\mod p&\t{if $3\mid y$,}
\\1\mod p&\t{if $3\nmid y$.}
\endcases\tag 3.3$$
Note that $2F_{n\pm 1}=L_n\pm F_n$. From Corollary 3.1 and the above
we deduce the result.

\pro{Theorem 3.5} Let $p$ be an odd  prime with $p\e 1,2,4,8\mod
{15}$.
\par $(\t{\rm i})$ If $p\e 1,4\mod{15}$ and so $p=x^2+15y^2$ with
$x,y\in\Bbb Z$, then
$$\sum_{k=0}^{[p/3]}\b{3k}k\f 1{3^k}\e \cases 1\mod p&\t{if
$3\mid y$,}
\\-(3x+5y)/(10y)\mod p&\t{if $3\mid y-x$.}
\endcases$$
\par $(\t{\rm ii})$ If $p\e 2,8\mod{15}$ and so $p=5x^2+3y^2$ with
$x,y\in\Bbb Z$, then
$$\sum_{k=0}^{[p/3]}\b{3k}k\f 1{3^k}\e \cases 1\mod p&\t{if
$3\mid y$,}
\\(3x-y)/(2y)\mod p&\t{if $3\mid y-x$.}
\endcases$$
\endpro
Proof. It is known that $U_n(3,1)=F_{2n}=F_nL_n$. Thus, putting
$a=b=\f 13$ in Theorem 3.3 we see that
$$\sum_{k=0}^{[p/3]}\f{\b{3k}k}{3^k}\e \cases -U_{\f{p-1}3-1}(3,1)
=-F_{\f{p-1}3-1}L_{\f{p-1}3-1}\mod p&\t{if $p\e 1\mod 3$,}
\\U_{\f{p+1}3+1}(3,1)
=F_{\f{p+1}3+1}L_{\f{p+1}3+1}\mod p&\t{if $p\e 2\mod 3$.}
\endcases$$
It is easily seen that
$$2F_{n\pm 1}=L_n\pm F_n\qtq{and} 2L_{n\pm 1}=5F_n\pm L_n.$$
Thus, if $p=x^2+15y^2\e 1,4\mod {15}$, using (3.2) we see that
$$F_{\f{p-1}3-1}L_{\f{p-1}3-1}=\f
14\big(L_{\f{p-1}3}-F_{\f{p-1}3}\big)\big(5F_{\f{p-1}3}-L_{\f{p-1}3}\big)
\e \cases -1\mod p&\t{if $3\mid y$,}
\\(3x+5y)/(10y)\mod p&\t{if $3\mid y-x$.}
\endcases$$
If $p=5x^2+3y^2\e 2,8\mod {15}$, using (3.3) we see that
$$F_{\f{p+1}3+1}L_{\f{p+1}3+1}=\f
14\big(L_{\f{p+1}3}+F_{\f{p+1}3}\big)\big(5F_{\f{p+1}3}+L_{\f{p+1}3}\big)
\e \cases 1\mod p&\t{if $3\mid y$,}
\\(3x-y)/(2y)\mod p&\t{if $3\mid y-x$.}
\endcases$$
Now combining all the above we obtain the result.

 \pro{Theorem 3.6} Let $p$ be an odd prime with $\sls p{13}=\sls p3$. Then
 $$\aligned&\sum_{k=0}^{[p/3]}\b{3k}k\f 1{(-3)^k}
 \\&\e\cases 1\mod p&\t{if $p=x^2+xy+88y^2$, $10x^2+7xy+10y^2$,}
 \\\qq\q&\t{or if $p=11x^2+xy+8y^2$,}
 \\ -(25x+10y)/(13y)\mod p&\t{if $p=25x^2+7xy+4y^2$,}
 \\(43x+12y)/(13y)\mod p&\t{if $p=43x^2+37xy+10y^2\not=43$,}
  \\-(5x+8y)/(13y)\mod p&\t{if $p=5x^2+3xy+18y^2\not=5$,}
 \\-(47x+9y)/(13y)\mod p&\t{if $p=47x^2+5xy+2y^2\not=47$.}
 \endcases\endaligned$$
 \endpro
Proof. Taking $b=\f 13$ and $a=-\f 13$ in Theorem 3.3 and applying
(2.2) we see that
$$\aligned&\sum_{k=0}^{[p/3]}\b{3k}k\f 1{(-3)^k}
 \\&\e\cases U_{\f{p-1}3-1}(3,-1)=-\f 12(3U_{\f{p-1}3}(3,-1)-V_{\f{p-1}3}(3,-1))
 \mod p&\t{if $\sls {13}p=1$,}
 \\-U_{\f{p+1}3+1}(3,-1)=-\f 12(3U_{\f{p+1}3}(3,-1)+V_{\f{p+1}3}(3,-1))
 \mod p&\t{if $\sls {13}p=-1$.}
 \endcases\endaligned$$
 Now applying [9, Corollary 6.7] we deduce the result.

 \pro{Theorem 3.7} Let $p$ be an odd prime with $\sls p3\sls p5\sls p{17}=1$. Then
 $$\aligned&\sum_{k=0}^{[p/3]}\b{3k}k(-3)^k
 \\&\e\cases 1\mod p&\t{if $p=x^2+xy+64y^2$, $3x^2+3xy+22y^2$,}
 \\\qq&\t{or if $p=8x^2+xy+8y^2$, $5x^2+5xy+14y^2$,}
 \\ -(171x+74y)/(85y)\mod p&\t{if $p=19x^2+7xy+4y^2\not=19$,}
 \\-(63x+65y)/(85y)\mod p&\t{if $p=7x^2+5xy+10y^2\not=7$,}
  \\-(63x+13y)/(17y)\mod p&\t{if $p=35x^2+5xy+2y^2$,}
 \\(99x-29y)/(17y)\mod p&\t{if $p=11x^2+3xy+6y^2\not=11$.}
 \endcases\endaligned$$
 \endpro
Proof. Taking $b=1$ and $a=-\f 13$ in Theorem 3.3 and applying (2.2)
 we see that
$$\aligned&\sum_{k=0}^{[p/3]}\b{3k}k(-3)^k
 \\&\e\cases U_{\f{p-1}3-1}(9,-1)=\f 12(V_{\f{p-1}3}(9,-1)-9
 U_{\f{p-1}3}(9,-1))\mod p&\t{if $\sls {85}p=1$,}
 \\-U_{\f{p+1}3+1}(9,-1)=-\f 12(V_{\f{p+1}3}(9,-1)+9
 U_{\f{p+1}3}(9,-1))\mod p&\t{if $\sls {85}p=-1$.}
 \endcases\endaligned$$
 Now applying [9, Corollary 6.9] we deduce the result.

\par Let $(u,v)$ be the greatest common divisor of integers $u$ and $v$.
 For $a,b,c\in\Bbb Z$ we use $[a,b,c]$ to denote the equivalence
class containing the form $ax^2+bxy+cy^2$. It is well known that
$$[a,b,c]=[c,-b,a]=[a,2ak+b,ak^2+bk+c]\qtq{for}k\in\Bbb Z.\tag 3.4$$
We also use $H(d)$ to denote the form class group of discriminant
$d$. Let $\o=(-1+\sqrt{-3})/2$. Following [6] and [9] we use $\sls
{a+b\o}m_3\ (3\nmid m)$ to denote the cubic Jacobi symbol. For a
prime $p>3$ and $k\in\Bbb Z_p$ with $k^2+3\not\e 0\mod p$, using [6,
Corollary 6.1] we can easily determine $\sls{k+1+2\o}p_3$. In
particular, by [6, Proposition 2.1] we have $\sls{1+2\o}p_3=1$.
\par For later convenience, following [9] we introduce the following
notation.
\pro{Definition 3.1} Suppose $u,v,d\inZ,\
dv(u^2-dv^2)\not=0$ and $(u,v)=1$.  Let $u^2-dv^2=2^{\a}3^rW(2\nmid
W, 3\nmid W)$ and let $w$ be the product of all distinct prime
divisors of $W$. Define
$$\aligned &k_2(u,v,d)=\cases 2&\t{if $d\e 2,3\mod 4$,}
\\2&\t{if $d\e 1\mod 8$, $\a>0$ and $\a\e 0,1\mod 3$,}
\\ 1&\t{otherwise,}\endcases
\\&k_3(u,v,d)=\cases 3^{\ord_3v+1}&\t{if $3\mid r$ and $3\nmid u$,}
\\ 9&\t{if $3\nmid r$ and $3\nmid u$,}
\\3&\t{if $3\nmid r-2$, $3\mid u$ and $9\nmid u$,}
\\1&\t{otherwise}
\endcases\endaligned$$
and $k(u,v,d)=k_2(u,v,d)k_3(u,v,d)w/(u,w).$
\endpro
 \pro{Lemma 3.3 ([9, Theorem 6.1 and Remark 6.1])} Let
$p>3$ be a prime, and $P,Q\inZ$ with $p\nmid Q$ and
$\ls{-3(P^2-4Q)}p=1$. Assume $P^2-4Q=df^2$ $(d,f\inZ)$ and
$p=ax^2+bxy+cy^2$ with $a,b,c,x,y\inZ$, $(a,6p\cdot 4Q/(P,f)^2)=1$
and $b^2-4ac=-3k^2d$, where $k=k(P/(P,f),f/(P,f),d)$. Then
$$U_{(p-\sls p3)/3}(P,Q)\e \cases 0\mod p&\t{if $\ls
{\f{bf}{(P,f)}-\f{kP}{(P,f)}(1+2\o)}{a}_3=1$,}
\\-\f{2ax+by}{kdfy}\ls {-Q}p(-Q)^{\f{p-\sls p3}6}\mod p&\t{if
$\ls{\f{bf}{(P,f)}-\f{kP}{(P,f)}(1+2\o)}{a}_3=\o$,}
\\\f{2ax+by}{kdfy}\ls {-Q}p(-Q)^{\f{p-\sls p3}
6}\mod p&\t{if
$\ls{\f{bf}{(P,f)}-\f{kP}{(P,f)}(1+2\o)}{a}_3=\o^2$}\endcases$$ and
$$V_{(p-\sls p3)/3}(P,Q)\e
\cases 2\ls p3\ls{-Q}p(-Q)^{\f{p-\sls p3}6}\mod p&\t{if
$\ls{\f{bf}{(P,f)}-\f{kP}{(P,f)}(1+2\o)}{a}_3=1$,}
\\ -\ls p3\ls{-Q}p(-Q)^{\f{p-\sls p3}6}\mod p
&\t{if $\ls{\f{bf}{(P,f)}-\f{kP}{(P,f)}(1+2\o)}{a}_3\not=1$.}
\endcases$$
Moreover, the criteria for $p\mid U_{(p-\sls p3)/3}(P,Q)$ and
$V_{(p-\sls p3)/3}(P,Q)\mod p$ are also true when $p=a$.
\endpro

 \pro{Theorem 3.8} Let $p>3$ be a prime with $\sls p{23}=1$. Then
 $$\sum_{k=0}^{[p/3]}\b{3k}k
 \e\cases 1\mod p&\t{if $p=x^2+xy+52y^2,\ 8x^2+7xy+8y^2$,}
  \\(39x-10y)/(23y)\mod p&\t{if $p=13x^2+xy+4y^2\not=13$,}
 \\-(87x+19y)/(23y)\mod p&\t{if $p=29x^2+5xy+2y^2\not=29$.}
 \endcases$$
 \endpro
Proof. Putting $a=b=1$ in Theorem 3.3 and applying $(2.2)$ we see
that
$$\sum_{k=0}^{[p/3]}\b{3k}k
 \e\cases (-3)^{\f {p-1}3+1}\cdot \f 16(9U_{\f{p-1}3}(9,3)-V_{\f{p-1}3}(9,3))\mod p
 &\t{if $3\mid p-1$,}
 \\-(-3)^{\f {p-2}3}\cdot \f 12(9U_{\f{p+1}3}(9,3)+V_{\f{p+1}3}(9,3))\mod p
 &\t{if $3\mid p-2$.}
 \endcases$$
 Since  $\sls p{23}=1$ we have $\sls{-3\cdot 69}p=\sls{-23}p=\sls
 p{23}=1$ and $\sls{69}p=\sls p3$. Thus $p$ is represented by some class in $H(-207)$. It is
 known that
 $$H(-207)=\{[1,1,52], [8,7,8],[4,1,13],[4,-1,13], [2,1,26],[2,-1,26]\}.$$
From (3.4) one can easily seen that $[2,-1,26]=[2,-5,29]=[29,5,2]$
and $[8,7,8]=[8,23,23]=[23,-23,8]$. Note that
$$\aligned
&\Ls{1-9(1+2\o)}1_3=1,\q\Ls{-23-9(1+2\o)}{23}_3=\Ls{1+2\o}{23}_3=1,
\\&\Ls{1-9(1+2\o)}{13}_3=\Ls{-3+1+2\o}{13}_3=\o,
\\&\Ls{5-9(1+2\o)}{29}_3=\ls{-7+1+2\o}{29}_3=\o^2.
\endaligned$$ Since $k(9,1,69)=1$ by Definition 3.1, putting $P=9$,
$Q=3$, $d=69$, $f=1$ and $k=1$ in Lemma 3.3 and applying the above
we see that
$$U_{\f{p-\sls p3}3}(9,3)\e\cases 0\mod p&\t{if $p=x^2+xy+52y^2$}
\\&\t{or $p=8x^2+7xy+8y^2$,}\\-\f{26x+y}{69y}(-3)^{\f{p-1}6}\mod p&\t{if
$p=13x^2+xy+4y^2\not=13$,}\\-\f{58x+5y}{69y}(-3)^{\f{p+1}6} \mod
p&\t{if $p=29x^2+5xy+2y^2\not=29$}
\endcases$$ and
$$V_{\f{p-\sls p3}3}(9,3)\e \cases 2(-3)^{(p-1)/6}\mod
p&\t{if $p=x^2+xy+52y^2$}
\\2(-3)^{(p+1)/6}\mod
p&\t{if $p=8x^2+7xy+8y^2$,}
\\-(-3)^{(p-1)/6}\mod p&\t{if
$p=13x^2+xy+4y^2$,}
\\-(-3)^{(p+1)/6}\mod p&\t{if
$p=29x^2+5xy+2y^2$.}
\endcases$$
Now combining all the above with the fact $(-3)^{(p-1)/2} \e \sls
{-3}p=\sls p3\mod p$ we deduce the result.

\pro{Theorem 3.9} Let $p>3$ be a prime with $\sls p{31}=1$. Then
 $$\sum_{k=0}^{[p/3]}\b{3k}k(-1)^k
 \e\cases 1\mod p&\t{if $p=x^2+xy+70y^2,9x^2+9xy+10y^2$}
 \\ &\q\t{or $p=8x^2+3xy+9y^2$,}
 \\(15x-14y)/(31y)\mod p&\t{if $p=5x^2+xy+14y^2\not=5$,}
 \\(21x-14y)/(31y)\mod p&\t{if $p=7x^2+xy+10y^2\not=7$,}
\\(57x-8y)/(31y)\mod p&\t{if $p=19x^2+5xy+4y^2\not=19$,}
 \\-(105x+17y)/(31y)\mod p&\t{if $p=35x^2+xy+2y^2$.}
 \endcases$$
 \endpro
Proof. Putting $a=-1$ and $b=1$ in Theorem 3.3 and applying $(2.2)$
 we see that
$$\sum_{k=0}^{[p/3]}\b{3k}k(-1)^k
 \e\cases 3^{\f {p-1}3+1}\cdot \f 1{(-6)}(9U_{\f{p-1}3}(9,-3)-V_{\f{p-1}3}(9,-3))\mod p
 &\t{if $3\mid p-1$,}
\\-3^{\f {p-2}3}\cdot \f 12(9U_{\f{p+1}3}(9,-3)+V_{\f{p+1}3}(9,-3))\mod p
 &\t{if $3\mid p-2$.}
 \endcases$$
 Since  $\sls p{31}=1$ we have $\sls{-3\cdot 93}p=\sls{-31}p=\sls
 p{31}=1$ and $\sls{93}p=\sls p3$. Thus $p$ is represented by some class in $H(-279)$. It is
 known that
 $$\align H(-279)&=\{[1,1,70], [9,9,10],[2,1,35],[2,-1,35],[5,1,14],[5,-1,14],
 \\&\qq[7,1,10],[7,-1,10],
[4,3,18],[4,-3,18],[8,3,9],[8,-3,9]\}.\endalign$$
  One can easily seen that $[2,-1,35]=[35,1,2]$, $[4,3,18]=[4,-5,19]=[19,5,4]$,
$[8,3,9]=[8,-29,35]=[35,29,8]$ and
$[9,9,10]=[10,-9,9]=[10,31,31]=[31,-31,10]$. Note that
$$\aligned
&\Ls{1-9(1+2\o)}1_3=1,\q\Ls{-31-9(1+2\o)}{31}_3=\Ls{1+2\o}{31}_3=1,
\\&\Ls{29-9(1+2\o)}{35}_3=\Ls{24+1+2\o}5_3\Ls{24+1+2\o}7_3=\o^2\cdot \o=1,
\\&\Ls{1-9(1+2\o)}5_3=\Ls{1+1+2\o}5_3=\o,
\q\Ls{1-9(1+2\o)}7_3=\ls{3+1+2\o}7_3=\o,
\\&\Ls{5-9(1+2\o)}{19}_3=\Ls{-9+1+2\o}{19}_3=\o,
\\&\Ls{1-9(1+2\o)}{35}_3=\Ls{-4+1+2\o}5_3\Ls{-4+1+2\o}7_3=\o\cdot \o=\o^2.
\endaligned$$ Since $k(9,1,93)=1$ by Definition 3.1, putting $P=9$,
$Q=-3$, $d=93$, $f=1$ and $k=1$ in Lemma 3.3 and applying the above
we see that
$$U_{\f{p-\sls p3}3}(9,-3)\e\cases 0\mod p&\t{if $p=x^2+xy+70y^2,9x^2+9xy+10y^2$}
\\&\t{or $p=8x^2+3xy+9y^2$,}
\\-\f{10x+y}{93y}\ls3p 3^{\f{p+1}6}\mod p&\t{if
$p=5x^2+xy+14y^2\not=5$,}
\\-\f{14x+y}{93y}\ls 3p 3^{\f{p-1}6} \mod
p&\t{if $p=7x^2+xy+10y^2\not=7$,}
\\-\f{38x+5y}{93y}\ls 3p 3^{\f{p-1}6}\mod p&\t{if
$p=19x^2+5xy+4y^2\not=19$,}
\\\f{70x+y}{93y}\ls 3p 3^{\f{p+1}6} \mod
p&\t{if $p=35x^2+xy+2y^2$}
\endcases$$ and
$$V_{\f{p-\sls p3}3}(9,-3)\e \cases
2\ls 3p 3^{(p-1)/6}\mod p&\t{if $p=x^2+xy+70y^2,9x^2+9xy+10y^2$,}
\\-2\ls 3p 3^{(p+1)/6}\mod
p&\t{if $p=8x^2+3xy+9y^2,$}
\\\ls 3p3^{(p+1)/6}\mod
p&\t{if $p=5x^2+xy+14y^2, 35x^2+xy+2y^2$,}
\\-\ls 3p 3^{(p-1)/6}\mod p&\t{if
$p=7x^2+xy+10y^2, 19x^2+5xy+4y^2$.}
\endcases$$
Now combining all the above  we deduce the result.
 \pro{Theorem
3.10} Let $p>3$ be a prime and $a\in\Bbb Z_p$ with
$\ls{a(4-27a)}p=-1$. Then $x\e \sum_{k=0}^{[p/3]}\b{3k}ka^k\mod p$
is the unique solution of the cubic congruence $(27a-4)x^3+3x+1\e
0\mod p$.
\endpro
Proof. As $\sls{a(4-27a)}p=-1$ we have
$\sls{81a^2-12a}p=\sls{-3}p\sls{a(4-27a)}p=-\sls{-3}p=-\sls p3$.
Thus putting $b=a$ in Theorem 3.3 we have
$$\sum_{k=0}^{[p/3]}\b{3k}ka^k\e \cases
-(3a)^{\f{p-1}3}U_{\f{p-1}3+1}(9a,3a)\mod p&\t{if $p\e 1\mod 3$,}
\\(-3a)^{\f{p-2}3+1}U_{\f{p+1}3-1}(9a,3a)\mod p&\t{if $p\e 2\mod 3$.}
\endcases$$
From [8, Theorem 2.1] or [9, Remark 6.1] we know that
$$U_{\f{p-\sls p3}3}(9a,3a)\e \f 1{27a-4}\Ls{-a}p(3a)^{\f{p-\sls
p3}6-1}(-3x^2+2x+18a)\mod p$$ and
$$V_{\f{p-\sls p3}3}(9a,3a)\e \Ls{3a}p(3a)^{\f{p-\sls
p3}6-1}(x^2-6a)\mod p,$$ where $x$ is the unique solution of the
congruence $X^3-9aX-27a^2\e 0\mod p$. Hence
$$\align&9aU_{\f{p-\sls p3}3}(9a,3a)+\Ls p3V_{\f{p-\sls p3}3}(9a,3a)
\\&\e \f 1{27a-4}\Ls{-a}p(3a)^{\f{p-\sls p3}6-1}
\big(9a(-3x^2+2x+18a)+(27a-4)(x^2-6a)\big)
\\&=-\f 2{27a-4}\Ls{-a}p(3a)^{\f{p-\sls p3}6-1}(2x^2-9ax-12a)
\mod p.\endalign$$ Now putting $b=a$ in Theorem 3.3 and applying
(2.2) and the above we deduce

$$\align \sum_{k=0}^{[p/3]}\b{3k}ka^k&\e -\Ls p3
(3a)^{\f{p-\sls p3}3}U_{\f{p-\sls p3}3+\sls p3}(9a,3a)\\&=-\Ls
p3(3a)^{\f{p-\sls p3}3}\cdot \f 1{2(3a)^{(1-\sls p3)/2}}
\Big(9aU_{\f{p-\sls p3}3}(9a,3a)+\Ls p3V_{\f{p-\sls
p3}3}(9a,3a)\Big)
\\&\e \Ls p3(3a)^{\f{p-\sls p3}3+\f{\sls p3-1}2}\f 1{27a-4}\Ls{-a}p
(3a)^{\f{p-\sls p3}6-1}(2x^2-9ax-12a)
\\&\e \f 1{3a(27a-4)}(2x^2-9ax-12a)\mod p.\endalign$$
As $x^3\e 9ax+27a^2\mod p$ we see that $$(2x^2-9ax-12a)(2x+9a)\e
3a(4-27a)x\mod p.$$ Hence
$$\sum_{k=0}^{[p/3]}\b{3k}ka^k\e
\f 1{3a(27a-4)}(2x^2-9ax-12a)\e -\f x{2x+9a}\mod p.$$ Set $x_0=-\f
x{2x+9a}$. Then $x=-\f{9ax_0}{2x_0+1}$ and
$$\align &(27a-4)x_0^3+3x_0+1\\&=(4-27a)\f{x^3}{(2x+9a)^3}-\f{3x}{2x+9a}+1
=-\f{27a(x^3-9ax-27a^2)}{(2x+9a)^3}\e 0\mod p.\endalign$$ So the
theorem is proved. \subheading{4. Congruences via combinatorial
sums}
\par \q For $m,n\in\Bbb N$ and $r\in\Bbb Z$ let
$$T_{r(m)}^n=\sum\Sb k\in\{0,1,\ldots,n\}\\ k\e r\mod m\endSb\b nk.$$
From [3, Corollary 1.8] or [11] we know that
$$T_{r(m)}^n=T_{n-r(m)}^n\qtq{and}
T_{r(m)}^{n+1}=T_{r(m)}^n+T_{r-1(m)}^n.\tag 4.1$$

\pro{Theorem 4.1} Let $p>3$ be a prime. Then
$$\sum_{k=0}^{[\f{p-1}6]}\f{\b{6k}{3k}}{(-64)^k}
\e\cases \f 13((-1)^{[\f{p+1}4]}+2(-1)^{\f{p-1}2})\mod p&\t{if
$3\mid p-1$,}
\\ \f
13((-1)^{[\f{p+1}4]}-(-1)^{\f{p-1}2})\mod p&\t{if $3\mid p-2$.}
\endcases$$
\endpro
Proof. It is known that (see for example ([1, (1.56)],[3, Theorem
1.1])
$$T_{0(3)}^n=\sum_{k=0}^{[n/3]}\b n{3k}=\cases \f 13(2^n+2(-1)^n)&\t{if $3\mid
n$,}\\\f 13(2^n-(-1)^n)&\t{if $3\nmid n$.}
\endcases\tag 4.2$$
Thus, taking $n=\f{p-1}2$ and using Lemma 2.3 we see that
$$\sum_{k=0}^{[\f{p-1}6]}\f{\b{6k}{3k}}{(-64)^k}
\e \sum_{k=0}^{[\f{p-1}6]}\b{\f{p-1}2}{3k} =\cases \f
13(2^{\f{p-1}2}+2(-1)^{\f{p-1}2})\mod p&\t{if $3\mid p-1$,}
\\ \f
13(2^{\f{p-1}2}-(-1)^{\f{p-1}2})\mod p&\t{if $3\mid p-2$.}
\endcases$$ To see the result, we note that $2^{\f{p-1}2}\e
(-1)^{[\f{p+1}4]}\mod p.$

\pro{Theorem 4.2} Let $p$ be an odd prime. Then
$$\sum_{k=0}^{[\f{p-1}8]}\f{\b{8k}{4k}}{4^{4k-1}}
\e\cases 1+(-1)^{\f{p-1}8}2^{\f{p+3}4}\mod p&\t{if $p\e 1\mod 8$,}
\\ -1+(-1)^{\f{p-3}8}2^{\f{p+1}4}\mod p&\t{if $p\e 3\mod 8$,}
\\-1\mod p&\t{if $p\e 5\mod 8$,}
\\ 1-(-1)^{\f{p-7}8}2^{\f{p+1}4}\mod p&\t{if $p\e 7\mod 8$.}
\endcases$$

\endpro
Proof. It is known that (see for example ([1, (1.58)], [3, Theorem
1.2])
$$T_{0(4)}^n=\sum_{k=0}^{[n/4]}\b n{4k}=\cases \f 12(2^{n-1}+(-1)^{[n/4]}2^{[n/2]})&\t{if $n\e 0,1
\mod 4$,}\\2^{n-2}&\t{if $n\e 2\mod 4$,}
\\\f 12(2^{n-1}-(-1)^{[n/4]}2^{[n/2]})&\t{if $n\e 3
\mod 4$.}
\endcases\tag 4.3$$
Thus, taking $n=\f{p-1}2$ and using Lemma 2.3 we see that
$$\sum_{k=0}^{[\f{p-1}8]}\f{\b{8k}{4k}}{(-4)^{4k}}
\e \sum_{k=0}^{[\f{p-1}8]}\b{\f{p-1}2}{4k} =\cases \f
12(2^{\f{p-1}2-1}+(-1)^{\f{p-1}8}2^{\f{p-1}4})\mod p&\t{if $8\mid
p-1$,}
\\ \f 12(2^{\f{p-1}2-1}+(-1)^{\f{p-3}8}2^{\f{p-3}4})\mod p&\t{if $8\mid p-3$,}
\\2^{\f{p-1}2-2}\mod p&\t{if $8\mid p-5$,}
\\ \f 12(2^{\f{p-1}2-1}-(-1)^{\f{p-7}8}2^{\f{p-3}4})\mod p&\t{if $8\mid p-7$.}
\endcases$$ To see the result, we
note that $2^{\f{p-1}2}\e (-1)^{[\f{p+1}4]}\mod p.$

\pro{Theorem 4.3} Let $p>3$ be a prime. Then
$$6\sum_{k=0}^{[\f{p-1}{12}]}\f{\b{12k}{6k}}{4^{6k}}
\e\cases 2\cdot 3^{(p-1)/4}+3\mod p&\t{if $p\e 1\mod{24}$,}
\\3^{(p-1)/4}-2\mod p&\t{if $p\e 5\mod{24}$,}
\\-2\mod p&\t{if $p\e 7\mod {24}$,}
\\ -3^{(p+1)/4}\mod p&\t{if $p\e 11\mod {24}$,}
\\1-2\cdot 3^{(p-1)/4}\mod p&\t{if $p\e 13\mod{24}$,}
\\-3^{(p-1)/4}\mod p&\t{if $p\e 17\mod{24}$,}
\\-3\mod p&\t{if $p\e 19\mod {24}$,}
\\ 2+3^{(p+1)/4}\mod p&\t{if $p\e 23\mod {24}$.}
\endcases$$

\endpro
Proof. From [3, Theorem 1.9] we know that
$$6\sum_{k=0}^{[n/6]}\b n{6k}-2^n=\cases
3^{(n+1)/2}+1&\t{if $n\e  \pm 1 \mod {12}$,}\\
-2&\t{if $n\e \pm 3\mod {12}$,}
\\-3^{(n+1)/2}+1&\t{if $n\e \pm 5
\mod {12}$,}
\\2(3^{n/2}+1)&\t{if $n\e 0\mod{12}$,}
\\3^{n/2}-1&\t{if $n\e \pm 2\mod{12}$,}
\\-3^{n/2}-1&\t{if $n\e \pm 4\mod{12}$,}
\\2(1-3^{n/2})&\t{if $n\e 6\mod{12}$.}
\endcases\tag 4.4$$
Thus, by the above and Lemma 2.3 we have
$$\aligned&6\sum_{k=0}^{[\f{p-1}{12}]}\f{\b{12k}{6k}}{(-4)^{6k}}-(-1)^{[\f{p+1}4]}
\\&\e 6\sum_{k=0}^{[\f{p-1}{12}]}\b{\f{p-1}2}{6k}-2^{\f{p-1}2}
 =\cases
 3^{(p+1)/4}+1\mod p&\t{if $p\e   23 \mod {24}$,}\\
-2\mod p&\t{if $p\e 7,19\mod {24}$,}
\\-3^{(p+1)/4}+1\mod p&\t{if $p\e  11
\mod {24}$,}
\\2(3^{(p-1)/4}+1)\mod p&\t{if $p\e 1\mod{24}$,}
\\3^{(p-1)/4}-1\mod p&\t{if $p\e 5\mod{24}$,}
\\-3^{(p-1)/4}-1\mod p&\t{if $p\e  17\mod{24}$,}
\\2(1-3^{(p-1)/4})\mod p&\t{if $p\e 13\mod{24}$.}
\endcases\endaligned$$ This yields the result.
\pro{Theorem 4.4} Let $p$ be a prime greater than $5$. Then
$$5\sum_{k=0}^{[\f{p-1}{10}]}\f{\b{10k}{5k}}{(-4)^{5k}}-(-1)^{[\f{p+1}4]}
\e\cases 4\cdot 5^{\f{p-1}4}\mod p&\t{if $p\e 1\mod{20}$,}
\\2\cdot 5^{\f{p+1}4}\mod p&\t{if $p\e
3\mod{20}$,}
\\5^{\f{p+1}4}\mod p&\t{if $p\e
7\mod{20}$,}
\\ 5^{\f{p-1}4}\mod p&\t{if $p\e
9\mod{20}$,}
\\ -5^{\f{p+1}4}\mod p&\t{if $p\e
11\mod{20}$,}
\\-2\cdot 5^{\f{p-1}4}\mod p&\t{if $p\e
13\mod{20}$,}
\\3\cdot 5^{\f{p-1}4}\mod p&\t{if $p\e
17\mod{20}$,}
\\-5^{\f{p+1}4}\mod p&\t{if $p\e
19\mod{20}$.}
\endcases$$
\endpro
Proof. Let
$$\Delta_5(r,n)=\cases 5T_{\f{n-1}2+r(5)}^n-2^n&\t{if $2\nmid n$,}
\\5T_{\f n2+r(5)}^n-2^n&\t{if $2\mid n$.}
\endcases$$
From [3, Theorem 1.6] we know that
$$\Delta_5(0,n)=2(-1)^nL_n,\ \Delta_5(\pm 1,n)=(-1)^nL_{n-1},\
\Delta_5(\pm 2,n)=(-1)^{n+1}L_{n+1}.$$ From this we deduce
$$5\sum_{5\mid k}\b{\f{p-1}2}k-2^{\f{p-1}2}=
\cases 2L_{\f{p-1}2}&\t{if $p\e 1\mod{20}$,}
\\-2L_{\f{p-1}2}&\t{if $p\e 3\mod{20}$,}
\\-L_{\f{p-3}2}&\t{if $p\e 7,19\mod{20}$,}
\\-L_{\f{p+1}2}&\t{if $p\e 9,13\mod{20}$,}
\\L_{\f{p+1}2}&\t{if $p\e 11\mod{20}$,}
\\L_{\f{p-3}2}&\t{if $p\e 17\mod{20}$.}
\endcases\tag 4.5$$
Using Lemma 2.3 we have
$$\sum_{k=0}^{[\f{p-1}{10}]}\f{\b{10k}{5k}}{(-4)^{5k}}\e
\sum_{k=0}^{[\f{p-1}{10}]}\binom{(p-1)/2}{5k}=\sum\Sb k=0\\ 5\mid
k\endSb^{(p-1)/2}\binom {(p-1)/2}k\mod p.$$ Note that
$L_{\f{p-3}2}=L_{\f{p+1}2}-L_{\f{p-1}2}$. In [11, Corollaries 1 and
2], the author and Z.W. Sun determined $L_{\f{p\pm 1}2}\mod p$.
 By the above and
[11, Corollaries 1 and 2] we deduce the result.

 \pro{Theorem 4.5}
Let $p\e 11\mod{20}$ be a prime. Then
$$\sum_{k=0}^{\f{p-11}{20}}\b{20k}{10k}\f 1{4^{10k}}\e
(-1)^{\f{p+1}4}\f 1{10}\mod p. $$
\endpro
Proof.  By Lemma 2.3 we have
$$10\sum_{k=0}^{\f{p-11}{20}}\b{20k}{10k}\f 1{4^{10k}}-(-1)^{\f{p+1}4}
\e
10\sum_{k=0}^{\f{p-11}{20}}\b{\f{p-1}2}{10k}-2^{\f{p-1}2}=10T_{0(10)}^{\f{p-1}2}
-2^{\f{p-1}2}\mod p.$$  According to [11, Theorem 1 and Corollary
1],
$$10T_{0(10)}^{\f{p-1}2}
-2^{\f{p-1}2}=-2L_{\f{p-1}2}\e 0\mod p.$$ Thus the result follows.

 \pro{Theorem 4.6} Let $p\e 13\mod{24}$ be a prime. Then
$$\sum_{k=0}^{(p-13)/24}\b{24k}{12k}\f
1{4^{12k}}\e \f 1{12}(1-2\cdot 3^{\f{p-1}4})\mod p$$ and
$$\sum_{k=0}^{(p-1)/12}\b{12k}{6k}\f {(-1)^k}{4^{6k}}\e 0\mod p.$$
\endpro
 Proof. Using Lemma 2.3 and the fact $2^{\f{p-1}2}\e
-1\mod p$ we have
$$12T_{0(12)}^{\f{p-1}2}-2^{\f{p-1}2}=12\sum_{k=0}^{[p/24]}
\b{\f{p-1}2}{12k}-2^{\f{p-1}2}\e 12\sum_{k=0}^{[p/24]}
\b{24k}{12k}\f 1{(-4)^{12k}}+1\mod p.$$ Since $\f{p-3}2\e
5\mod{12}$, by (4.1) and [12, Theorem 2] we have
$$12T_{0(12)}^{\f{p-3}2}-2^{\f{p-3}2}=1-3^{\f{p-1}4}+(-1)^{\f{p-5}8}(2^{\f{p-1}4}-
V_{\f{p-5}4}(4,1))$$ and
$$12T_{-1(12)}^{\f{p-3}2}-2^{\f{p-3}2}=12T_{\f{p-1}2(12)}^{\f{p-3}2}-2^{\f{p-3}2}
=1-3^{\f{p-1}4}-(-1)^{\f{p-5}8}(2^{\f{p-1}4}- V_{\f{p-5}4}(4,1)).$$
Thus, using (4.1) we obtain
$$12T_{0(12)}^{\f{p-1}2}-2^{\f{p-1}2}=12T_{0(12)}^{\f{p-3}2}-2^{\f{p-3}2}
+12T_{-1(12)}^{\f{p-3}2}-2^{\f{p-3}2}=2(1-3^{\f{p-1}4}).$$ Hence
$$\sum_{k=0}^{[p/24]}\b{24k}{12k}\f 1{4^{12k}}\e \f
1{12}(-1+2(1-3^{\f{p-1}4}))
\mod p.$$
 Since
$$6\sum_{k=0}^{(p-1)/12}\b{12k}{6k}\f 1{4^{6k}}+6\sum_{k=0}^{(p-1)/12}\b{12k}{6k}\f {(-1)^k}{4^{6k}}
=12\sum_{k=0}^{(p-13)/24}\b{24k}{12k}\f 1{4^{12k}},$$ by the above
and Theorem 4.3 we deduce the remaining result. 
\enddocument
\bye